\documentclass[11 pt]{article}
 \usepackage{amssymb}
\usepackage{verbatim}
\usepackage{latexsym}
\usepackage{eucal}
\def\proof{\mbox {\bf Proof.}\quad}
\def\vfi{\varphi}
\newtheorem{sub}{\name}[section]
\newcommand{\bs}{
\begin{sub}}
\newcommand{\es}{
\end{sub}}
\newcommand{\bsl}[1]{
\begin{sub}\label{#1}}
\newcommand{\bth}[1]{\def\name{Theorem}
\begin{sub}\label{t:#1}}
\newcommand{\blemma}[1]{\def\name{Lemma}
\begin{sub}\label{l:#1}}
\newcommand{\bcor}[1]{\def\name{Corollary}
\begin{sub}\label{c:#1}}
\newcommand{\bdef}[1]{\def\name{Definition}
\begin{sub}\label{d:#1}}
\newcommand{\bprop}[1]{\def\name{Proposition}
\begin{sub}\label{p:#1}}
\newcommand{\brem}[1]{\def\name{Remark}
\begin{sub}\label{r:#1}}

\newcommand{\rth}[1]{Theorem~\ref{t:#1}}
\newcommand{\rlemma}[1]{Lemma~\ref{l:#1}}

\newcommand{\rprop}[1]{Proposition~\ref{p:#1}}


\newcommand{\BA}{
\begin{array}}
\newcommand{\EA}{
\end{array}}
\newcommand{\BAN}{\renewcommand{\arraystretch}{1.2}
\setlength{\arraycolsep}{2pt}
\begin{array}}
\newcommand{\BAV}[2]{\renewcommand{\arraystretch}{#1}
\setlength{\arraycolsep}{#2}
\begin{array}}
\newcommand{\BSA}{
\begin{subarray}}
\newcommand{\ESA}{
\end{subarray}}
\newcommand{\BAL}{
\begin{aligned}}
\newcommand{\EAL}{
\end{aligned}}
\newcommand{\BALG}{
\begin{alignat}}
\newcommand{\EALG}{
\end{alignat}}%% the abbrev. does not work with latex2e
\newcommand{\BALGN}{
\begin{alignat*}}
\newcommand{\EALGN}{
\end{alignat*}}%% the abbrev. does not work with latex2e
\newcommand{\note}[1]{\textit{#1.}\hspace{2mm}}
\newcommand{\Proof}{\note{Proof}}
\newcommand{\qeda}{\hspace{10mm}\hfill $\square$}
\newcommand{\qed}{\\${}$ \hfill $\square$}
\newcommand{\Remark}{\note{Remark}}
\newcommand{\abs}[1]{\left |#1\right |}%% adjustable vertical 
\def\angb<#1>{\langle #1 \rangle}%% angle bracket
\newcommand{\opname}[1]{\mbox{\rm #1}\,}

\newcommand{\dist}{\opname{dist}}
\newcommand{\myfrac}[2]{{\displaystyle \frac{#1}{#2} }}
\newcommand{\myint}[2]{{\displaystyle \int_{#1}^{#2}}}

\newcommand{\prt}{
\partial}

\def\ga{\alpha}     \def\gb{\beta}       \def\gg{\gamma}
       \def\gd{\delta}      \def\ge{\epsilon}
\def\gth{\theta}                         
\def\gf{\phi}           
            \def\gl{\lambda}
\def\gm{\mu}        \def\gn{\nu}         \def\gp{\pi}
    \def\gr{\rho}        
\def \gs{\sigma}       
      \def\gw{\omega}
                
\def\Gg{\Gamma}     \def\Gd{\Delta}

      \def\CN{{\mathcal N}}
      
\def\CA{{\mathcal A}}

   \def\BBN {\mathbb N}    
   \def\BBR {\mathbb R}

\catcode`\@=11

\def\la{\lambda}
\def\de{\delta}
\def\vep{\varepsilon}

\def\eqalign#1{\null\,\vcenter{\openup1\jot \m@th
  \ialign{\strut\hfil$\displaystyle{##}$&$\displaystyle{{}##}$\hfil
     &&\strut$\displaystyle{##}$\hfil&$\displaystyle{{}##}$
     \hfil\crcr#1\crcr}}\,}
\newcommand{\rife}[1]{(\ref{#1})}
\def\qed{{\unskip\nobreak\hfil\penalty50
          \hskip2em\hbox{}\nobreak\hfil\mbox{\rule{1ex}{1ex} \qquad}
   \parfillskip=0pt
   \finalhyphendemerits=0\par\medskip}}

\begin{document}
\title {\bf Separable p-harmonic functions in a cone and related quasilinear equations on manifolds}
\author{{\bf\large Alessio Porretta}\\
 {\it Dipartimento di Matematica}\\ {\it Universita di Roma Tor Vergata, Roma}\\[2mm]
{\bf\large Laurent V\'eron}\\
{\it Laboratoire de Math\'ematiques et Physique Th\'eorique}\\ {\it Universit\'e Fran\c cois Rabelais, Tours}}

\date{}
\maketitle

{\small{\bf Abstract} In considering a class of quasilinear elliptic equations on a Riemannian manifold with nonnegative Ricci curvature, we give a new proof of Tolksdorf's result on the construction of separable $p$-harmonic functions in a cone. }
\vspace{1mm}
\hspace{.05in}
%\parbox{4.5in} {{\small }}

\noindent {\it \footnotesize 1991 Mathematics Subject 
Classification}. {\scriptsize 35K60
}.\\
{\it \footnotesize Key words}. {\scriptsize  
}
%\vspace{1cm}
\section {Introduction}
\setcounter{equation}{0} 
Let $(r, \gs)$ be the spherical coordinates in $\BBR^N$. If $u$ is a harmonic function in $\BBR^N\setminus\{0\}$ written under the separable form 
\begin{equation}\label{s1}
u(x)=r^{-\gb}\gw( \gs)
\end{equation}
 it is straightforward to check that $\gw$ is an eigenfunction of the Laplace-Beltrami operator $-\Gd_{_{S^{N-1}}}$ on the unit sphere $S^{N-1}\subset\BBR^N$ and $\gb$ is a root of
\begin{equation}\label{h1}
X^2-(N-2)X-\gl=0,
\end{equation}
where $\gl\geq 0$ is the corresponding eigenvalue. The function $\gw$ is called a spherical harmonic and its properties  are well--known, since such functions are  the restrictions to the sphere of homogeneous harmonic polynomials. More generally, if $C_{S}\subset\BBR^N$ is the cone with vertex $0$ and opening $S\subsetneq S^{N-1}$, there exist positive harmonic functions $u$ in $C_{S}$ under the form (\ref{s1}) which vanish on $\prt C_{S}\setminus\{0\}$ if and only if $\gb$ is a root of (\ref{h1}), where, in that case, $\gl:=\gl_{_{S}}$ is the first eigenvalue of $-\Gd_{_{S^{N-1}}}$ in $W^{1,2}_{0}(S)$. These separable harmonic functions play a fundamental role in the discription of isolated interior or boundary singularities of solutions of second order linear elliptic equations. If the Laplace equation is replaced by the $p$-Laplace equation 
\begin{equation}\label{p-L}
-\Gd_{p}u:=-{\rm div}\left(|D u|^{p-2}D u\right)=0,
\end{equation}
($p>1$), the same question of existence of separable p-harmonic functions, i.e. solutions of (\ref{p-L}) in the form \rife{s1}, was considered by Krol \cite{Kr}, Tolksdorf \cite{To},  Kichenassamy and V\'eron
 \cite{KV}. If $u$ in \rife{s1} is p-harmonic, then the function $\gw$ must be a solution of the spherical p-harmonic equation, 
 \begin{equation}\label{s-p-h}
-div\left((\gb^2\gw^2+|\nabla' \gw|^2)^{p/2-1}\nabla' \gw\right)=\gb(\gb (p-1)+p-N)(\gb^2\gw^2+|\nabla' \gw|^2)^{p/2-1}\gw,
\end{equation}
on $S^{N-1}$, where $\nabla'$ and $div$ are respectively the covariant derivative identified with the \lq\lq tangential gradient\rq\rq\  and the divergence operator acting on vector fields on $S^{N-1}$.  Two special cases arise when either $p=2$ or $N=2$: if $p=2$, (\ref{s-p-h}) is just an eigenvalue problem
 \begin{equation}\label{s-2-h}
-\Gd'\gw=\gb(\gb +2-N)\gw,
\end{equation}
where $\Gd'$ is the Laplace-Beltrami operator on $S^{N-1}$. When  $N=2$,  equation (\ref{s-p-h}) becomes
 \begin{equation}\label{s-p-h1}
-\left((\gb^2\gw^2+| \gw_{\gth}|^2)^{p/2-1}\gw_{\gth}\right)_{\gth}=\gb(\gb (p-1)+p-2)(\gb^2\gw^2+|\gw_{\gth}|^2)^{p/2-1}\gw,
\end{equation}
where $\gth\in [0,\gp]$. Introducing the new unknown $\gf:=\gw_{\gth}/\gw$, (\ref{s-p-h1}) is transformed into a separable equation,
 \begin{equation}\label{s-p-h2}
-\left((\gb^2+\gf^2)^{p/2-1}\gf\right)_{\gth}=\left((p-1)\gf^2+\gb(\gb(p-1)+p-2)\right)(\gb^2+\gf^2)^{p/2-1}.
\end{equation}
This equation was completely integrated by Krol \cite{Kr} in the case $\gb<0$, and Kichenassamy and V\'eron
 \cite{KV} in the case $\gb>0$. It turns out that for any integer $k>0$ there exist two couples $(\tilde \gb_{k},\tilde \gf_{k})$ and $( \gb_{k},\gf_{k})$ where $\tilde \gb_{k}<0$, $\gb_{k}>0$, and $\tilde \gf_{k}$ and $\gf_{k}$ are anti-periodic solutions of the corresponding equation  (\ref{s-p-h2}). Furthermore $\tilde \gf_{k}$ and $\gf_{k}$ are uniquely determined, up to an homothety. 
 
 A remarkable breakthrough was realized by Tolksdorf \cite {To} when he proved that for any smooth domain $S\subset S^{N-1}$ there exists a couple $(\gb,\gf)$ where $\gb<0$ and $\gf\in C^{1}(\bar S)$ is positive in $S$, vanishes on $\prt S$ and solves (\ref{s-p-h}) in $S$. Furthermore $\gb:=\tilde \gb_{_{S}}$ is unique and $\gf$ is determined up to a multiplicative constant. Tolksdorf's result is obtained by constructing a p-harmonic function $u$ in the cone $C_{_{S}}$ generated by $S$ with a compactly supported boundary data and by proving, thanks to a kind of Harnack inequality up to the boundary, the \lq\lq equivalence principle\rq\rq, that the asymptotic behaviour of $u$ is self-similar. Later on  the existence of a couple $(\gb,\gf)$, with $\gb:= \gb_{_{S}}>0$ and $\gf$, as above, positive solution of (\ref{s-p-h}) in $S$ vanishing on $\prt S$  is proved by the same method in \cite {Ve1}, therefore we shall refer to the two cases $\gb>0$ and $\gb<0$ as Tolksdorf's results. The structure of these spherical $p$-harmonic functions is studied in \cite{BoV} when $p=N$. These regular ($\gb<0$) and singular ($\gb>0$) separable $p$-harmonic functions play a fundamental role in describing the behaviour of solutions of quasilinear equations near  a regular or singular boundary point \cite{Kr},\cite{KM},\cite{BBoV1},\cite{BoV2}.\\
 
 In this article, we give a new proof of Tolksdorf's results, entirely different from his. Actually, performing a change of variable, we embed our problem into a much wider class of quasilinear equations. 
Indeed, 
if  $\gw\in W^{1,p}_{0}(S)$ is a positive solution of (\ref{s-p-h}) in $S\subset S^{N-1}$, which vanishes on $\prt S$, then 
%it is $C^1$ in $\bar S$. Furthermore the normal derivative (relative to the Riemannian structure) $u_{\gn}$ at the boundary is negative. because the equation is never degenerate the standard regularity theory yields to $v\in C{\infty}(S)\cap C^2(\bar S)$. 
the function $v$ defined by
$$v=-\myfrac{1}{\gb}\ln\gw
$$
solves
 \begin{equation}\label{mainl}\left\{ \BA {l}
 -div\left(\left(1+|\nabla' v|^2\right)^{p/2-1}\nabla' v\right)+\gb(p-1)\left(1+|\nabla' v|^2\right)^{p/2-1}|\nabla' v|^2\\[2mm]\phantom{---.\left(\left(1+|\nabla' v|^2\right)^{p/2-1}\right)}
 =-\left(\gb(p-1)+p-N\right)\left(1+|\nabla' v|^2\right)^{p/2-1}\quad\mbox{ in }S
 \\[2mm]\phantom{\left(\left(\right)^{p=x,2}\nabla' v\right)}
 \lim_{ \gs\to\prt S}v( \gs)=\infty.
 \EA\right.
 \end {equation}
Notice that this equation is never degenerate and $v$  is $C^2$ (actually  $C^{\infty}$) in $S$ and satisfies the equation and the boundary condition in classical sense. 
%As a developed form,
 %\begin{equation}\label{mainI}\left\{ \BA {l}
 %-\Gd'v-(p-2)\myfrac{D^2v \nabla'v.\nabla'v}{1+|\nabla'v|^2}+\gb(p-1)|\nabla'v|^2=-\left(\gb(p-1)+p-N\right)\quad\mbox{ in }S
% \\[2mm]\phantom{\left(\left(---------\right)^{p=x,}\nabla' v\right)}
 %\lim_{ \gs\to\prt S}v( \gs)=\infty.
 %\EA\right.
 %\end {equation}
 Our construction of solutions of (\ref{s-p-h}) relies on  a careful study 
 of the quasilinear problem  \rife{mainl}, and on the interpretation of the constant in the right hand side of \rife{mainl} as  an  \lq\lq ergodic constant\rq\rq. Furthermore, having an intrinsic independent interest,  this study  will be  performed on any compact smooth subdomain of a Riemannian manifold, without refering to the p-Laplace equation \rife{p-L}. Our main result is the following:\\ 
 
 \noindent{\bf Theorem A}. {\it Let $(M,g)$ be a $d$-dimensional Riemannian manifold with nonnegative Ricci curvature, and let $\nabla$ and $div_{g}$ be respectively the covariant derivative and the divergence operator on $M$. Then for any compact smooth subdomain $S\subset M$ and any $\gb>0$ there exists a unique positive constant $\gl_{\gb}$ such that the problem 
 \begin{equation}\label{main1}\left\{ \BA {l}
 -div_{g}\left(\left(1+|\nabla v|^2\right)^{p/2-1}\nabla v\right)+\gb(p-1)\left(1+|\nabla v|^2\right)^{p/2-1}|\nabla v|^2\\[2mm]\phantom{-div_{g}\left(\left(1+|\nabla v|^2\right)^{p/2-1}\nabla v\right)}
 =-\gl_{\gb}\left(1+|\nabla v|^2\right)^{p/2-1}\quad\quad\quad\mbox{ in }\;S
 \\[2mm]\phantom{-div_{g}\left(--;\nabla v\right)}
 \lim_{x\to\prt S}v(x)=\infty.
 \EA\right.
 \end {equation}
 admits a solution $v\in C^2( S)$. Furthermore, $v$ is unique up to an additive constant.}\\

\noindent  The result of Theorem A is the typical statement of an  ergodic problem, indeed  the constant $\la_\beta$ can be seen as  the unique ergodic constant for the equation obtained after  dividing by $\left(1+|\nabla v|^2\right)^{p/2-1}$ (see  \rife{gen1}). Observe also that \rife{main1} may be reformulated if we set $\gw=e^{-\gb v}$, then $\gw$ is a solution of
 \begin{equation}\label{s-p-h'}\left\{ \BA {l}
-div_{g}\left((\gb^2\gw^2+|\nabla\gw|^2)^{p/2-1}\nabla \gw\right)=\gb\gl_{\gb}(\gb^2\gw^2+|\nabla \gw|^2)^{p/2-1}\gw\quad\mbox{ in }\; S\\[2mm]
\phantom{-div_{g}\left((\gb^2\gw^2+|\nabla\gw|^2)^{p/2-1}\nabla\right)}
\gw=0\quad\mbox{ on }\;\prt S
 \EA\right.
 \end {equation}
When $p=2$, problem (\ref{s-p-h'}) reduces to an eigenvalue problem since $\gb\gl_{\gb}=\gl_{_{1}}(S)$, the principal eigenvalue of the Laplace--Beltrami operator in $S$. In that case the connection between \rife{main1} and \rife{s-p-h'} dates back to the stochastic interpretation of principal eigenvalues
(see  e.g. \cite{Ho}, \cite{LL}). 
 In the nonlinear framework with $p\neq 2$, by  proving  that the mapping $\gb\mapsto\gl_{\gb}$ is continuous, decreasing and tends to $\infty$ as $\gb\to 0^+$, we conclude that the equation $ \gl_{\gb}=(\gb (p-1)+p-d-1)$ has  a unique positive solution. As a consequence we generalize Tolksdorf's result as follows.\\
 
 \noindent{\bf Theorem B}. {\it Under the assumptions of Theorem A, for any compact smooth subdomain $S$ of $M$ 
 there exists a unique $\gb:=\gb_{_{S}}>0$ such that the problem
\begin{equation}\label{O1}\left\{\BA{l}
-div_{g}\left((\gb^2\gw^2+|\nabla \gw|^2)^{p/2-1}\nabla \gw\right)=\gb\left(\gb(p-1)+p-d-1\right)(\gb^2\gw^2+|\nabla \gw|^2)^{p/2-1}\gw\,\,\mbox{in }S\\[2mm]
\phantom{;;-\left((\gb^2\gw_{0}^2+|\nabla \gw_{0}|^2)^{p/2-1}\nabla \gw_{0}\right)}
\gw=0\quad\mbox{on }\prt S,
\EA\right.
\end{equation}
admits a positive solution $\gw\in C^1(\bar S)\cap C^2(S)$. Furthermore $\gw$ is unique up to an homothethy.}

 Of course, we obtain similarly that  for $\gb<0$ there exists a unique $\gb:=\tilde\gb_{_{ S}}<0$ such that $ \gl_{\gb}= (\gb (p-1)+p-d-1)$. Tolksdorf's results then follow as a particular case by taking $(M,g)=(S^{N-1},g_{0})$, where $S^{N-1}$ is equipped with the standard metric $g_{0}$ induced by the Euclidean structure in $\BBR^N$.
 
 %A starting remark, which justifies the intrinsic point of view is that, when $p=2$, problem (\ref{s-p-h}) reduces to the eigenvalue problem (\ref{s-2-h}) on the sphere since $\gb_{S}\gl_{\gb_{S}}=\tilde\gb_{S}\gl_{\tilde\gb_{S}}=\gl_{_{1}}(S)$. Actually Tolksdorf called (\ref{s-p-h}) a nonlinear eigenvalue problem. Our proof  consists precisely  in adopting the strategy to  look, for a given $\gb$, at the coefficient in the right hand side of (\ref{main1}) as to an eigenvalue.  In order to  study  its dependence on $\beta$, we go back to a possible constructive way of finding eigenvalues, which originally dates back to the stochastic interpretation of principal eigenvalues  and is connected with ergodic problems and boundary blow--up solutions  (see  e.g. \cite{Ho}, \cite{LL}). In this sense, our proof of Theorem A contains related results on ergodic problems which may have an independent  interest.\\

%%%%%%%%%%%%%%%%%%%%%%%%%%
%%%%%%%%%%%%%%%%%%%%%%%%%%
\section{The singular case}
\setcounter{equation}{0}

%and  $S^{N-1}$ is equiped with the standard metric $g_{0}$ induced by the Euclidean structure in $\BBR^N$. If  $\gw\in W^{1,p}_{0}(S)$ is a positive solution of (\ref{s-p-h}) in $S\subset S^{N-1}$, which vanishes on $\prt S$, then it is $C^1$ in $\bar S$. Furthermore the normal derivative (relative to the Riemannian structure) $u_{\gn}$ at the boundary is negative. because the equation is never degenerate the standard regularity theory yields to $v\in C{\infty}(S)\cap C^2(\bar S)$. 
In the following, we consider a general geometric setting and we recall some elements of Riemannian geometry (see e.g.  \cite{doc}, \cite{pole}). Let $(M,g)$ be a complete $d$-dimensional Riemannian manifold with metric tensor $g=(g_{ij})$, inverse $g^{-1}=(g^{ij})$ and determinant $|g|$. 
If $X$ and $Y$ are two tangent vector fields to $M$, we denote by
$$X.Y=\sum_{{ij}}{}g_{ij}(x)X^iY^j
$$
their scalar product in the tangent space $T_{x}M$. 
Let $x_j$, $j=1,...,d$, be a local system of coordinates: if  $u\in C^1(M)$, the gradient of $u$, quoted by $\nabla u$,  is the vector field with 
%covariant components $\nabla_{i}u=u_{x_{i}}$ and contravariant 
components
$\left(\nabla u\right)^{i}=\sum_{{k}}{}g^{ik}u_{x_{k}}$. Therefore
$$\nabla u.\nabla u=|\nabla u|^2=\sum_{{ij}}{}g^{ij}(x)u_{x_{i}}u_{x_{j}}.
$$
If $X=(X^{i})$ is a $C^1$ vector field on $M$, the divergence of $X$ is defined by
$$div_{g}X=\myfrac{1}{\sqrt{|g|}}\sum_{{k}}{}\left(\sqrt{|g|}X^k\right)_{x_{k}}.
$$
Recalling that, in local coordinates, the Christoffel symbols are
$$
\Gg_{ij}^{k}=\myfrac{1}{2}\sum_{{l}}{}
\left(\myfrac{\prt g_{jl}}{\prt x_i}+\myfrac{\prt g_{li}}{\prt x_j}
-\myfrac{\prt g_{ij}}{\prt x_l}\right)g^{lk},
$$
the second covariant derivatives of a $C^2$ function $u$ are
$$\nabla_{ij}u=u_{x_{i}x_{j}}-\sum_{{k}}{}\Gg_{ij}^{k}u_{x_k},
$$
while the Hessian is the $2$-tensor $D^2u=(\nabla_{ij}u)$. Finally, $\Gd_{g}u=trace( D^2u)=div_{g}\nabla u$ is the Laplace-Beltrami operator on $M$, locally expressed by
$$\Gd_{g}u=\myfrac{1}{\sqrt{|g|}}\sum_{{ij}}
\myfrac{\prt}{\prt x_{i}}\left(\sqrt{|g|}\,g^{ij}\myfrac{\prt u}{\prt x_{j}}\right)
=\sum_{{ij}}\myfrac{\prt}{\prt x_{i}}
\left(g^{ij}\myfrac{\prt u}{\prt x_{j}}\right)+ \sum_{{ijk}}\Gg^k_{ik}\,g^{ij}
\myfrac{\prt u}{\prt x_{j}}.$$
We denote by $Ricc_{g}$ the Ricci curvature tensor of the metric $g$.  In particular, if $(M,g)=(S^{N-1},g_{0})$,  then $Ricc_{g_{0}}=(N-1)g_{0}$.

 In all the sequel $p>1$ is a real number. We prove next the result of Theorem A, which we restate  here for the reader's convenience. 

\bth{dev-th} Let $S\subset M$ be a smooth bounded open domain of $M$ such that $Ricc_{g}\geq 0$ on $S$. Then for any $\gb>0$ there exists a unique $\gl_{\gb}>0$ such that there exists a function $v\in C^2(S)$ satisfying
 \begin{equation}\label{gen1}\left\{ \BA {l}
 -\Gd_{g}v-(p-2)\myfrac{D^2v \nabla v.\nabla v}{1+|\nabla v|^2}+\gb(p-1)|\nabla v|^2=-\gl_{\gb}\quad\mbox{ in }S
 \\[2mm]\phantom{\left(\left(---------\right)^{p=x,}\nabla v\right)}
 \lim_{x\to\prt S}v( x)=\infty.
 \EA\right.
 \end {equation}
Furthermore, $v$ is unique up to an additive constant.
\es
\Proof  As in the usual approach to ergodic problems, we start  by considering the  problem
 \begin{equation}\label{gen2}\left\{ \BA {l}
 -\Gd_{g}v_{\ge}-(p-2)\myfrac{D^2v_{\ge} \nabla v_{\ge}.\nabla v_{\ge}}{1+|\nabla v_{\ge}|^2}+\gb(p-1)|\nabla v_{\ge}|^2+\ge v_{\ge}=0\quad\mbox{ in }S
 \\[2mm]\phantom{\left(\left(------+\ge v_{\ge}----\right)^{p=x,}\nabla v\right)}
 \lim_{x\to\prt S}v_{\ge}( x)=\infty,
 \EA\right.
 \end {equation}
 where $\ge>0$, and then we study  the limit when $\ge\to 0$.\medskip
 
\noindent{\it Step 1: Construction of super and sub solutions. } Since $\prt S$ is $C^2$, the   distance function $\gr(x)=\dist (x,\prt S)$, where the distance is the geodesic distance, is a positive $C^2$ function is some relative neighborhood $\CN_{\gd}=\{x\in M:|\dot\gr(x)|<\gd\}$ of $\prt S$; here $\dot\gr(x)$ is the signed distance, equal to $\pm\gr(x)$ according $x\in S$ or $x\in M\setminus S$. Then $|\nabla\dot\gr(x)|=1$ in $\CN_\gd$. We extend $\dot\gr$ outside $\CN_\gd$ into a $C^2(M)$ function $\tilde\gr$. Next we consider the function
\begin{equation}\label{barr}
\bar u(x)=-\myfrac{1}{\gb}\ln(\tilde\gr(x))-M_{0}\tilde\gr(x)+\myfrac{M_{1}}{\ge}\qquad\forall x\in  S, 
\end{equation}
where the $M_{j}>0$ are to be chosen later on. Then
$$\nabla \bar u(x)=-\myfrac{1}{\gb\tilde\gr(x)}\left(1+\gb M_{0}\tilde\gr(x)\right)\nabla\tilde\gr(x),
$$
$$|\nabla \bar u(x)|^2=\myfrac{1}{\gb^2\tilde\gr^2(x)}\left(1+\gb M_{0}\tilde\gr(x)\right)^2\, |\nabla\tilde\gr(x)|^2.
$$
Notice that this last identity implies
$$|\nabla \bar u(x)|^2=\myfrac{1+2\gb M_{0}\gr(x)+O(\gr^2(x))}{\gb^2\gr^2(x)}\quad\mbox{as }\gr(x)\to 0.
$$
Next
$$\BA {l}
-\Gd_{g}\bar u-(p-2)\myfrac{D^2\bar u\nabla \bar u.\nabla \bar u}{1+|\nabla \bar u|^2}=-\Gd_{g}\bar u-\myfrac{(p-2)}{2}\myfrac{\nabla(|\nabla \bar u|^2).\nabla \bar u}{1+|\nabla \bar u|^2}\\[4mm]
%%%%%%
\phantom{\myfrac{D^2\bar u\nabla \bar u.\nabla \bar u}{1+|\nabla \bar u|^2}}
=-\myfrac{|\nabla \tilde\gr|^2}{\beta\,\tilde\gr^2}
+ \myfrac{\Delta_{g} \tilde\gr}{\beta\, \tilde\gr}(1+\beta\,M_0 \tilde\gr) 
-\myfrac{(p-2)}{\beta\,\tilde\gr^2}\myfrac{|\nabla \tilde\gr|^4(1+\beta \,M_0\,\tilde \gr)^3}{\beta^2\,\tilde\gr^2+|\nabla \tilde\gr|^2(1+\beta\,M_0\,\tilde\gr)^2}\\[4mm]
\phantom{\myfrac{D^2\bar u\nabla \bar u.\nabla \bar u}{1+|\nabla \bar u|^2}}
%%%%
+\myfrac{(p-2)}{2\beta\,\tilde\gr}\myfrac{(1+\beta\,M_0\,\tilde\gr)\,\nabla \tilde\gr\cdot\nabla\left[ |\nabla \tilde\gr|^2(1+\beta\,M_0\,\tilde\gr)^2\right]}{\beta^2\,\tilde\gr^2+|\nabla \tilde\gr|^2(1+\beta\,M_0\,\tilde\gr)^2}.
\EA$$
After some lengthy but standard computations, one obtains the following relation
\begin{equation}\label{barr1}\BA {l}
-\Gd_{g}\bar u-(p-2)\myfrac{D^2\bar u\nabla \bar u.\nabla \bar u}{1+|\nabla \bar u|^2}+\gb(p-1)|\nabla\bar u|^2+\ge\bar u\\[2mm]
\phantom{---------}
=\myfrac{1}{ \tilde\gr}\left(\myfrac{ \Gd_{g}\tilde\gr}{\gb}-\myfrac\vep\beta \tilde\gr\ln(\tilde \gr)+2(p-1)M_{0}|\nabla\tilde\gr|^2\right) +\psi_{\gb}(x)+M_{1},
\EA\end{equation}
where $\psi_{\gb}$ is a  function depending on $\gb$ (and on $M_0$), but 
which remains bounded on $S$, uniformly when $\gb$ remains in a compact subset of $(0,\infty)$. Since $|\nabla \tilde\gr|=1$ near the boundary, it is possible to choose  $M_{0}$ and $M_{1}$ such that $\bar u$ defined by (\ref{barr}) is a supersolution for (\ref{gen2}). Moreover, $M_{0}$ and $M_{1}$  can be chosen independent of  $\gb$ whenever it varies on a compact subset of $(0,\infty)$.  

One finds similarly that the function
\begin{equation}\label{under}
\underline u(x)=-\myfrac{1}{\gb}\ln(\tilde\gr(x))+M_{0}\tilde\gr(x)-\myfrac{M_{1}}{\ge}\qquad\forall x\in S, 
\end{equation}
is a subsolution of (\ref{gen2}), with $M_{0}$ and $M_{1}$ chosen as for $\bar u$. Moreover, for $0<h<\gd$, we can approximate $\bar u$ and $\underline u$ respectively from above and from below by
\begin{equation}\label{bar-h}
\bar u_{h}(x)=-\myfrac{1}{\gb}\ln(\tilde\gr(x)-h)-M_{0}(\tilde\gr(x)-h)+\myfrac{M_{1,h}}{\ge},
\end{equation}
\begin{equation}\label{under-h}
\underline u_{h}(x)=-\myfrac{1}{\gb}\ln(\tilde\gr(x)+h)+M_{0}(\tilde\gr(x)+h)-\myfrac{M_{1,h}}{\ge}, 
\end{equation}
which are, respectively, a  supersolution in $\{x\in S:\gr(x)>h\}$ and  a subsolution in $S$. Together with the comparison principle, these super and sub solutions will be used to derive estimates on the solutions of (\ref{gen2}).\medskip
 
\noindent{\it Step 2: Basic estimates. } In this part, by using the classical Bernstein's method (\cite{ber}), we derive the fundamental gradient estimate for the solutions $u\in C^2(S)$ of
 \begin{equation}\label{gen3}
 -\Gd_{g}u-(p-2)\myfrac{D^2u \nabla u.\nabla u}{1+|\nabla u|^2}+\gb(p-1)|\nabla u|^2+\ge u=0\quad\mbox{ in }S.
\end{equation}
We recall the Weitzenb\"ock formula (see e.g. \cite{Berger}):
\begin{equation}\label{wei}
\myfrac{1}{2}\Gd_{g}|\nabla u|^2=|D^2u|^2+\nabla(\Gd_{g}u).\nabla u+Ricc_{g}(\nabla u,\nabla u),
\end{equation}
and the Cauchy-Schwarz inequality for $D^2u$
$$|D^2u|^2\geq\myfrac{1}{d}|\Gd_{g}u|^2.
$$
Let $m=\inf\{Ricc_{g}(\nabla u,\nabla u):|\nabla u|=1\}\geq 0$, then
\begin{equation}\label{wei1}\myfrac{1}{2}\Gd_{g}|\nabla u|^2\geq 
\myfrac{1}{d}|\Gd_{g}u|^2+m|\nabla u|^2+\nabla(\Gd_{g}u).\nabla u.
\end{equation}
If we set $z=|\nabla u|^2$, we can re-write (\ref{gen3}) as
 \begin{equation}\label{gen4}
 \Gd_{g}u=-\myfrac{(p-2)}{2}\myfrac{\nabla z.\nabla u}{1+|\nabla u|^2}+\gb(p-1)z+\ge u\qquad\mbox{ in }S.
\end{equation}
Using the fact that
$$\nabla(\nabla z.\nabla u).\nabla u=D^2z\nabla u.\nabla u+\myfrac{1}{2}|\nabla z|^2,
$$
we obtain
$$\BA {l}
\nabla(\Gd_{g}u).\nabla u=-\myfrac{(p-2)}2\,\myfrac{D^2z\nabla u.\nabla u}{1+|\nabla u|^2}-\myfrac{(p-2)}{4}\myfrac{|\nabla z|^2}{1+|\nabla u|^2}
+\myfrac{(p-2)}2\myfrac{(\nabla z.\nabla u)^2}{(1+|\nabla u|^2)^2}\\
\phantom{2\nabla(\Gd_{g}u).\nabla u=-(p-2)\myfrac{D^2z\nabla u.\nabla u}{1+z}-\myfrac{(p-2)}{2}\myfrac{|\nabla z|^2}{1+|\nabla u|^2}}+\gb(p-1)\nabla z.\nabla u+\ge z.
\EA$$
Since, from (\ref{gen4})
$$|\Gd_{g}u|^2\geq c_{0}z^2-c_{1}\left((\ge u^-)^2+\myfrac{(\nabla z.\nabla u)^2}{(1+|\nabla u|^2)^2}\right),
$$
we derive from (\ref{wei1})
$$\BA {l}
\Gd_{g}z+(p-2)\myfrac{D^2z\nabla u.\nabla u}{1+|\nabla u|^2}\geq
\myfrac{2c_{0}z^2}{d}-\myfrac{2c_{1}}{d}\left((\ge u^-)^2 +\myfrac{(\nabla z.\nabla u)^2}{(1+|\nabla u|^2)^2}\right)+2(m+\ge)z\\
\noalign{\medskip}
\phantom{\Gd_{g}z+(p-2)\myfrac{D^2z\nabla u.\nabla u}{1+|\nabla u|^2}}-\myfrac{(p-2)}{2}\myfrac{|\nabla z|^2}{1+|\nabla u|^2}+(p-2)\myfrac{(\nabla z.\nabla u)^2}{(1+|\nabla u|^2)^2}+2\gb(p-1)\nabla z.\nabla u,
\EA$$
which yields, by Young's inequality and the fact that $z=|\nabla u|^2$,
\begin{equation}\label{gen5}
-\Gd_{g}z-(p-2)\myfrac{D^2z\nabla u.\nabla u}{1+|\nabla u|^2}+C_{0}z^2+2(m+\ge)z\leq C_{1}\myfrac{|\nabla z|^2}{1+z}+C_{2}
\end{equation}
for some positive constants $C_{j}$ ($j=0,1,2$), eventually depending on $\beta$, with  the   constant $C_2$ also depending on $\|\ge u^-\|_\infty$. Next we introduce the operator $\CA$ defined by
\begin{equation}\label{gen6}
\CA(z)=-\Gd_{g}z-(p-2)\myfrac{D^2z\nabla u.\nabla u}{1+|\nabla u|^2},
\end{equation}
which can be written, in local coordinates, as
\begin{equation}\label{gen6'}\BA {l}
\displaystyle\CA(z)=-\sum_{ij}g^{ij}z_{x_{i}x_{j}}-(p-2)\sum_{ij}\myfrac{\displaystyle\sum_{kl}g^{ik}u_{x_k}g^{jl}u_{x_l}}{1+\displaystyle\sum_{kl}g^{kl}u_{x_k}u_{x_l}}z_{x_{i}x_j}\\
\phantom{\CA(z)}\displaystyle
-\sum_{ijk}\left(\left(\Gg^k_{ik}g^{ij}+g^{ij}_{x_i}\right)z_{x_j}
+(p-2)\Gg^k_{ij}\myfrac {\displaystyle\sum_{lm}g^{im}u_{x_m}g^{jl}u_{x_l}}{1+\displaystyle\sum_{lm}g^{ml}u_{x_m}u_{x_l}}z_{x_k}\right)\\
\phantom{\CA(z)}\displaystyle
=-\sum_{ij}a_{ij}z_{x_{i}x_{j}}+ \sum_{i}b_{i}z_{x_{i}},
\EA\end{equation}
where the $a_{ij}$ are uniformly elliptic and bounded and the $b_{i}$ are bounded: indeed, it holds
$$
\min(p-1,1)g^{ij}\xi_i\xi_j\leq a_{ij}\xi_i\xi_j\leq \max(1,p-1)g^{ij}\xi_i\xi_j\,. 
$$ 
%Furthermore, if $X$ is any nonnegative matrix $X$ and every vector $\xi$ we have
%$$\min(1,p-1)tr (X) \leq tr (X)+(p-2)\myfrac{X\xi.\xi}{1+|\xi|^2}
%\leq \max(1,p-1)tr (X) .
%$$
Therefore from \rife{gen5} $z$ is a positive subsolution of an equation of the type
\begin{equation}\label{gen7}
\CA(z)+h(z)+g(z)|\nabla z|^2=f,
\end{equation}
where $g(z)=- C_1(1+z)^{-1}$,  $h(z)=2(m+\ge )z+C_{0}z^2$ and $f=C_2$. Since $m\geq 0$, $g$ and $h$ are increasing functions of the nonnegative variable $z$,   it follows that the comparison principle holds between super and sub-solutions of 
\begin{equation}\label{gen8}
-\Gd_{g}z-(p-2)\myfrac{D^2 z\nabla u.\nabla u}{1+|\nabla u|^2}+C_{0}z^2+2(m+\ge )z-C_{1}\myfrac{|\nabla z|^2}{1+z}=C_{2}.
\end{equation}
Standard computations show that, if  $\gl$ and $\gm$ are positive constants large enough, the function
$$\bar z(x)=\myfrac{\gl}{{\tilde \gr}^2 (x)}+\gm
$$
is a supersolution of  \rife{gen8}, which in addition blows up on $\prt S$.  We conclude that any bounded subsolution of (\ref{gen8}) satisfies $z(x)\leq \bar z(x)$, and therefore any subsolution by replacing $S$ by $\{x\in S:\gr(x)>h\}$ and $\tilde \gr(x)$ by $\tilde \gr(x)-h$. 

Finally, we proved that any $u\in C^2(S)$ which is solution of (\ref{gen3}) satisfies
\begin{equation}\label{gen9}
|\nabla u(x)|\leq \myfrac{L_0}{\tilde \gr(x)}+L_1\quad\forall x\in S,
\end{equation}
for some constants $L_0$, $L_1$ depending on $\|\vep\,u^-\|_\infty$.   Moreover, $L_0$ and $L_1$ can be chosen   uniformly bounded with respect to $\gb$, provided $\gb$ remains in a compact subset of $(0,\infty)$.

To conclude with the estimates on solutions of (\ref{gen3}), it is classical from the theory of quasilinear elliptic equations (see e.g. \cite{GT}) that local Lipschitz estimates imply local $C^{2,\ga}$ estimates since the equation is smooth  and uniformly elliptic.\medskip
 
\noindent{\it Step 3: Existence for the approximate equation. } As in \cite{LL}, we consider, for $n\in\BBN$ the solution $v_{n,\ge}:=v$ of
 \begin{equation}\label{gen2n}\left\{ \BA {l}
 -\Gd_{g}v-(p-2)\myfrac{D^2v \nabla v.\nabla v}{1+|\nabla v|^2}+\gb(p-1)|\nabla v|^2+\ge v=0\quad\mbox{ in }S
 \\[2mm]\phantom{\left(\left(------;;--+\ge v_{\ge}----\right)^{p=x,}\nabla v\right)}
v( x)=n\quad\mbox{ on }\prt S,
 \EA\right.
 \end {equation}
By previous steps,  the following estimates hold in $S$.
 \begin{equation}\label{es1}
0\leq v_{n,\ge}(x)\leq-\myfrac{1}{\gb}\ln\tilde \gr(x)-M_{0}\tilde \gr(x)+\myfrac{M_{1}}{\ge},
\end{equation} 
\begin{equation}\label{es2}
|\nabla v_{n,\ge}(x)|\leq\myfrac{L_0}{\tilde \gr(x)}+L_1.
\end{equation}
Moreover the sequence $\{v_{n,\ge}\}$ is bounded in $C^{2,\ga}_{loc}(S)$, which ensures the local compactness of the gradients. Since $n\mapsto v_{n,\ge}$ is increasing, there exists $v_{\ge}=\lim_{n\to\infty}v_{n,\ge}$  and $v_\ge$  is a solution of (\ref{gen2}) which satisfies (\ref{es1}) and  (\ref{es2}).\medskip
 
\noindent{\it Step 4: The ergodic limit. }  From Step 1, by  comparison with $\bar u_{h}$ and $\underline u_{h}$ defined in \rife{bar-h}--\rife{under-h} (and letting $h\to0$), we know that their holds in $S$:
 \begin{equation}\label{es4}
-\myfrac{1}{\gb}\ln\tilde\gr(x)+M_{0}\tilde\gr(x)-\myfrac{M_{1}}{\ge}
\leq v_{\ge}(x)\leq -\myfrac{1}{\gb}\ln\tilde\gr(x)-M_{0}\tilde \gr(x)+\myfrac{M_{1}}{\ge}.
\end{equation}
Therefore $\ge v_{\ge}$ is locally bounded in $S$. Since $\nabla v_{\ge}$ is locally bounded too in $S$, $\ge_{n} v_{\ge_{n}}$ converges to some constant $\gl_{0}\geq 0$ for some sequence $\{\ge_{n}\}$ in the $C_{loc}$-topology of $S$. We fix $x_{0}\in S$ and set $w_{\ge}:=v_{\ge}(x)-v_{\ge}(x_{0})$. Because $w_{\ge}$ is locally bounded in $C^1_{loc}(S)$ and $w_{\ge}$ satisfies 
 \begin{equation}\label{gen10}
 -\Gd_{g}w_{\ge}-(p-2)\myfrac{D^2w_{\ge} \nabla w_{\ge}.\nabla w_{\ge}}{1+|\nabla w_{\ge}|^2}+\gb(p-1)|\nabla w_{\ge}|^2+\ge w_{\ge}=-\ge v_{\ge}(x_{0})\quad\mbox{ in }S
\end{equation}
the regularity theory for elliptic equations implies that 
$w_{\ge}$ is locally bounded in $C^{2,\ga}(S)$. Up to an extraction of subsequence, there exists $w_{0}=\lim_{n\to\infty}w_{\ge_{n}}$, and $w_{0}$ is a solution of
 \begin{equation}\label{gen11}
 -\Gd_{g}w_{0}-(p-2)\myfrac{D^2w_{0} \nabla w_{0}.\nabla w_{0}}{1+|\nabla w_{0}|^2}+\gb(p-1)|\nabla w_{0}|^2=-\gl_{0}\quad\mbox{ in }S.
\end{equation}
The only question which remains to be proved is that $w_{0}$ blows-up at the boundary. We set
$$\underline\psi (x)=-\myfrac{1}{\gb}\ln\tilde \gr(x)+M_{0}\tilde \gr(x),
$$
and get, with same computations as in (\ref{barr1}),
\begin{equation}\label{barrx}\BA {l}
-\Gd_{g}\underline\psi-(p-2)\myfrac{D^2\underline\psi\nabla \underline\psi.\nabla \underline\psi}{1+|\nabla \underline\psi |^2}+\gb(p-1)|\nabla\underline\psi|^2+\ge\underline\psi\\[4mm]
\phantom{---------}
=
\myfrac{1}{ \tilde\gr}\left(\myfrac{ \Gd_{g}\tilde\gr}{\gb}-\myfrac\vep\beta \tilde\gr\ln(\tilde \gr) -2(p-1)M_{0}|\nabla\tilde\gr|^2\right)+\psi_{\gb}(x),
\EA\end{equation}
where $\psi_{\gb}$ is  a bounded function (depending on $\beta$, $M_0$). Noticing that $|\nabla\tilde\gr|=1$ in a neighborhood of $\prt S$, and that $ \ge v_{\ge}(x_{0})$ is uniformly bounded, we can choose $M_{0}$, $\gr_0$ such that the function $\underline\psi$ is a subsolution of (\ref{gen10}) in $\{x\in S\,:\,0<\gr(x)<\gr_0\}$. Since, whenever $\gr(x)=\gr_0$, we have  $w_{\ge}(x)\geq -c_0$ for some $c_0>0$  (due to the gradient estimate for $v_\ge$), and since $\underline\psi-c$ is still a subsolution for any positive constant $c$, we derive
\begin{equation}\label{barrw}
w_{\ge}(x)\geq -\myfrac{1}{\gb}\ln\tilde \gr(x)+M_{0}\tilde \gr(x)-c\qquad\forall x\mbox{ s.t. }\gr(x)\leq\gr_0.
\end{equation}
Letting $\ge$ tend to $0$ implies that $\lim_{x\to\prt S}w_{0}(x)=\infty$.\medskip
 
\noindent{\it Step 5: Uniqueness of the ergodic limit. } 
We claim that there exists a unique constant $\gl_{0}> 0$ such that there exists $v_{0}\in C^2(S)$ solution of
\begin{equation}\label{gen12}\left\{ \BA {l}
 -\Gd_{g}v_{0}-(p-2)\myfrac{D^2v_{0} \nabla v_{0}.\nabla v_{0}}{1+|\nabla v_{0}|^2}+\gb(p-1)|\nabla v_{0}|^2=-\gl_{0}\quad\mbox{ in }S
 \\[2mm]\phantom{\left(\left(----+v_{0} ----\right)^{p=x,}\nabla v_{0}\right)}
 \lim_{x\to\prt S}v_{0}( x)=\infty.
 \EA\right.
 \end {equation}
To this purpose, it will be useful the following 
\vskip0.3em

\blemma{regu}
A   function $v_0\in C^2(S)$ is solution of \rife{gen12} if and only if   the function $\gw_{0}=e^{-\gb v_{0}}\in C^2(S)\cap C(\bar S)$ is a  solution of
 \begin{equation}\label{clo3}\left\{\BA {l}
-div_{g}\left((\gb^2\gw_{0}^2+|\nabla \gw_{0}|^2)^{p/2-1}\nabla \gw_{0}\right)=\gb\gl_{0}(\gb^2\gw_{0}^2+|\nabla \gw_{0}|^2)^{p/2-1}\gw_{0}\quad\mbox{in }S\\[2mm]
\phantom{;;-\left((\gb^2\gw_{0}^2+|\nabla \gw_{0}|^2)^{p/2-1}\nabla \gw_{0}\right)}
\gw_{0}=0\quad\mbox{on }\prt S.
\EA\right.\end{equation}
Moreover, $\gw_0\in C^{1,\gg}(\bar S)$ for some $\gg>0$, and $\prt_{\gn}\gw_{0}<0$ on $\prt S$.
\es

\proof
Let $v_0\in C^2(S)$ be a solution of \rife{gen12}. As in the previous steps, considering the functions 
$$\underline\phi(x)=-\myfrac{1}{\gb}\ln\tilde\gr(x)+M_{0}\tilde\gr(x)-M^{*}\quad\mbox{and }\;\bar\phi(x)=-\myfrac{1}{\gb}\ln\tilde\gr(x)-M_{0}\tilde\gr(x)+M^{*},
$$
which appear to be respectively a sub and a super-solution for (\ref{gen12}) in $\{x\,:\, \gr(x)<\gd\}$ for some $\gd>0$ small enough (where $M^{*}$ depends on the value of 
$v_0$ on the set $\{x\in S:\gr(x)=\gd\}$), we obtain
\begin{equation}\label{clo}
\abs{v_{0}(x)+\myfrac{\ln\tilde \gr(x)}{\gb}}\leq M^{*}.
\end{equation}
By the gradient estimates of Step 2, there holds
\begin{equation}\label{clo2}
|\nabla v_0(x)|\leq \myfrac{L_0}{\tilde \gr(x)}+L_1\,.
\end{equation}
Now set $\gw_{0}=e^{-\gb v_{0}}$, then $\gw_{0}\in W^{1,\infty}(S)\cap C(\bar S)$ solves the problem \rife{clo3}. 
By the  regularity theory for  degenerate equations of $p$--Laplacian type (see the Appendix, \rth{liebe} and related references), we can deduce   that  $\gw_{0}\in C^{1,\gg}(\bar S)$. Moreover, since \rife{clo} implies 
\begin{equation}\label{barrier}
e^{-\beta M^*}\leq \frac{\gw_0}{\rho(x)}\leq e^{\beta M^*}
\end{equation}
we deduce that $\prt_{\gn}\gw_{0}<-e^{-\beta M^*}<0$ on $\prt S$. As a consequence, since $\gw_0\in C^{1}(\bar S)$ and is positive in $S$, we deduce that problem \rife{clo3} is uniformly elliptic, so that the classical regularity theory applies to give $\gw_0\in  C^{2,\ga}(S)$.

Of course, the converse is also true: given a solution $\gw_0$ of \rife{clo3}, clearly $v_0=-\frac1\beta \ln \gw_0$ is a  solution of \rife{gen12}. 
\qed

Assume now that there exist two ergodic constants, $\gl_{1}$ and $\la_2$, associated with two  solutions $v_{1}$, $v_2$, and let   correspondingly  $\gw_{i}=e^{-\gb v_{i}}$ be solutions of \rife{clo3}. Notice that multiplying (\ref{clo3}) by 
$\gw_{0}$ and integrating on $S$, we get  actually $\gl_{0}>0$.  Thus $\gl_{i}>0$ and, say, $\gl_{2}>\gl_{1}$. 

Since $\gw_{1}/\gw_{2}\in L^\infty(S)$ (from estimate \rife{barrier}), we denote
$$
\gth=\sup_{S}\myfrac{\gw_{1}}{\gw_{2}}.
$$
Because equation (\ref{clo3}) is homogeneous we can assume that $\gth=1$ and either there exists $x_{0}\in S$ such that $\gw_{1}(x_{0})=\gw_{2}(x_{0})$, $\nabla\gw_{1}(x_{0})=\nabla\gw_{2}(x_{0})$ and $\gw_{1}(x)\leq \gw_{2}(x)$ for $x\in \bar S$, or $\gw_{1}(x)<\gw_{2}(x)$ for  $x\in S$ and there exists $x_{0}\in\prt S$ such that  $\prt_{\gn}\gw_{1}(x_{0})=\prt_{\gn}\gw_{2}(x_{0})$. In the first case, it turns out that the function $z=v_{1}-v_{2}$ is nonnegative in
$S$, achieves a minimum at $x_{0}\in S$ and satisfies
$$-\Gd_{g}z(x_{0})-(p-2)\myfrac{D^2z(x_{0})\nabla v_{1}(x_{0}).\nabla v_{1}(x_{0})}{1+|\nabla v_{1}(x_{0})|^2}=\gl_{2}-\gl_{1}>0,
$$
which is impossible because of ellipticity. In the second case,  we have $\prt _{\gn}(\gw_1-\gw_2)(x_0)=0$, whereas $\gw_1-\gw_2$ is negative in $S$ and $(\gw_1-\gw_2)(x_0)=0$. Since the problem \rife{clo3} is uniformly elliptic (recall that the functions $\gw_i$ satisfy $(\beta^2 w_i^2+ |\nabla \gw_i|^2)>0$ on $\overline S$) this contradicts Hopf maximum principle. Therefore $\gw_{1}=\gw_{2}$, which implies $\gl_{1}=\gl_{2}$ by the equation. Thus the ergodic constant is unique. 

In a similar way one can prove  that $\omega_0$ is unique up to a multiplicative constant, and so $v_0$ is unique up to an additive constant (as a consequence, the whole sequence
$w_{\ge}$, constructed in Step 4, converges to $ w_{0}$ as $\ge\to 0$).

However, the uniqueness of $v_0$ can be proved with  a more general argument, concerning directly problem \rife{gen12}, which is  a variant as well as a generalization of previous uniqueness results for explosive solutions. Since it can have its own interest, we present it here.

First of all, we recall  that any $C^2$ function   $v_0$ solution of \rife{gen12} satisfies \rife{clo} and \rife{clo2}. Moreover, by \rlemma{regu} we have that $\gw_0=e^{-\beta v_0}\in  C^{1}(\bar S)$ and $\prt_{\gn}\gw_{0}<0$ on $\prt S$, hence, using that $\nabla v_0=- \frac{e^{\beta v_0}}\beta \nabla \gw_0$ and the estimate  \rife{clo} we conclude that there exists a constant $\sigma>0$ such that, in a neighborhood of $\partial S$
\begin{equation}\label{sotto}
|\nabla v_0| \geq \frac{\sigma}{\rho(x)}\,.
\end{equation}
In addition, it is possible 
to deduce from \rife{clo}--\rife{clo2} that there exists a constant $C_0>0$ such that  
\begin{equation}\label{2der}
|D^2 v_0|\leq \frac{C_0}{\tilde\rho^2(x)} \qquad \forall x\in S\,.
\end{equation} 
Indeed, take $x_0\in S$ and let $\rho_0=\frac{\rho(x_0)}2$, where we recall that $\rho(x_0)={\rm dist}(x_0,\partial S)$. Then consider (in a local neighborhood of $x_0$) the rescaled function 
$$
u_0(\xi)=v_0(x_0+\gr_0\,\xi)+\frac{\ln \gr_0}\beta\,,
$$
for $\xi \in B(0,1)$.  Note that $\rho(x_0+\gr_0\,\xi)\in (\rho_0, 3\rho_0)$ so  that \rife{clo2} and \rife{sotto} imply $\frac\sigma3 \leq |Du_0|\leq L_0+L_1\,\rho_0$.
Since $v_0$ is a solution of \rife{gen12}, a simple scaling in the local coordinates gives that $u_0$ is a solution of
$$
 -\Gd_{g}u_{0}-(p-2)\myfrac{D^2u_{0} \nabla u_{0}.\nabla u_{0}}{\rho_0^2+|\nabla u_{0}|^2}+\gb(p-1)|\nabla u_{0}|^2=-\gl_{0}\, \gr_0^2\quad\mbox{for $\xi \in B(0,1)$}
$$
with a  slight abuse of notation since now, in local coordinates, the derivatives are taken with respect to the variable $\xi$. Since the second order operator is uniformly elliptic (as in \rife{gen6}--\rife{gen6'}), by the classical regularity theory (e.g. see \cite{GT}, Theorem 13.6 to deduce the H\"older estimates for $Du_0$ and then apply the Schauder estimates, Chapter 6) we have that 
$$
|D^2 u_0(\xi)|\leq C\qquad\forall\xi\in B\left(0,{\frac 12}\right)
$$ 
where $C$ is a  constant depending on $\sup_{B(0,1)} \left(|u_0|+ |Du_0|\right)$. Using the estimates \rife{clo}--\rife{clo2} we can bound this last quantity only depending on $M^*$, $L_0$, $L_1$, hence we conclude that $|D^2 u_0(0)|\leq C$, which gives \rife{2der}. 

Now, take two solutions $v_1$, $v_2$ of \rife{gen12} corresponding to $\la_1$, $\la_2$ with, say,  $\la_1\leq \la_2$. We adapt now an argument in \cite{LL}: consider the 
function $\hat v= \theta v_2$, for $\theta<1$, and compute
$$
\left\{ \BA {l}
 -\Gd_{g} \hat v-(p-2)\myfrac{D^2 \hat v \nabla \hat v.\nabla \hat v}{1+|\nabla \hat v|^2}+\gb(p-1)|\nabla \hat v|^2=-\theta\gl_{2} \\[2mm]
\qquad + (1-\theta^2)\theta(p-2)\frac{D^2v_2\nabla v_2.\nabla v_2}{(1+|\nabla v_2|^2)(1+\theta^2|\nabla v_2|^2)}
- (1-\theta) \theta\beta(p-1)|\nabla v_2|^2 \EA\right.
$$ 
Using \rife{clo2}, \rife{2der} and \rife{sotto}, we know that 
$$
\left |\frac{D^2v_2\nabla v_2.\nabla v_2}{(1+|\nabla v_2|^2)(1+\theta^2|\nabla v_2|^2)}\right| \leq C\qquad \hbox{in $S$}
$$
hence $\hat v$ satisfies
$$
-\Gd_{g} \hat v-(p-2)\myfrac{D^2 \hat v \nabla \hat v.\nabla \hat v}{1+|\nabla \hat v|^2}+\gb(p-1)|\nabla \hat v|^2\leq   -\gl_{1} -(1-\theta)[\beta(p-1)\theta |\nabla v_2|^2- C (\theta+\theta^2)(p-2)-\la_1 ].
$$
Thanks to \rife{sotto}, we deduce that there exists $\de>0$, independent on $\theta$, such that $\hat v$ satisfies
$$
-\Gd_{g} \hat v-(p-2)\myfrac{D^2 \hat v \nabla \hat v.\nabla \hat v}{1+|\nabla \hat v|^2}+\gb(p-1)|\nabla \hat v|^2\leq - \gl_{1}
$$
 in $\{x\in S\,:\rho(x)<\de\}$. However, from the estimate \rife{clo} which holds for $v_1$ and $v_2$ we have that $v_1-\hat v\to +\infty$ as $\rho(x)\to 0$, hence $v_1-\hat v$ has a minimum in $\{x\in S\,:\rho(x)<\de\}$ and, by standard maximum principle, it is reached when $\rho(x)=\de$.
Letting $\theta\to 1$, we conclude that 
 $$
 \min\{(v_1-v_2)(x),\,\, x\,:\, \rho(x)\leq \de\}= \min\{(v_1-v_2)(x),\,\, x\,:\, \rho(x)=\de\}\,.
 $$
 On the other hand, looking at the equations of $v_1$, $v_2$ in $\{x\in S\,:\, \rho(x)>\de\}$, we also know (again by maximum principle) that 
 $$
 \min\{(v_1-v_2)(x),\,\, x\,:\, \rho(x)\geq \de\}= \min\{(v_1-v_2)(x),\,\, x\,:\, \rho(x)=\de\}
 $$ 
 hence $v_1-v_2$ should have a global minimum reached at  a point $x_0\in S$ such that $\rho(x_0)=\de$. Since $x_0$ lies inside the domain, and  the function $z=v_1-v_2$  satisfies a smooth elliptic equation around $x_0$, using the strong maximum principle  we conclude that $v_1- v_2$ is constant. This proves  the uniqueness, up to a constant,  of the  solution of \rife{gen12}, and at the same time also the  uniqueness of the ergodic constant ($\la_1=\la_2$, as we already proved before).
\qeda\medskip

\brem{}  \rm The argument  used in the last step of the previous proof also provides  a general uniqueness result for  explosive solutions of 
\begin{equation}\label{expl}
\left\{ \BA {l}
 -\Gd_{g}v-(p-2)\myfrac{D^2v \nabla v.\nabla v}{1+|\nabla v|^2}+\gb(p-1)|\nabla v|^2+ \ge v=f\quad\mbox{ in }S
 \\[2mm]\phantom{\left(\left(---------\right)^{p=x,}\nabla v\right)}
 \lim_{x\to\prt S}v( x)=\infty\,.
 \EA\right.
 \end {equation}
Precisely, if  $f$ is  a Lipschitz function, and   $\ge>0$,  the problem \rife{expl} has a unique solution $v\in C^2(S)$.  To our knowledge, such a  result is new even in the euclidean setting $M=\BBR^N$.
\es

\vskip1em
We proceed now studying how the ergodic constant  $\la_\beta$ depends on $\beta$, which will lead to the proof of Theorem B.
\vskip1em

\bprop {map} Under the assumptions of \rth{dev-th}, the mapping $\gb\mapsto\gl_{\gb}$ is continuous and  decreasing from $(0,\infty)$ in $(0,\infty)$, and it verifies
\begin{equation}\label{lim}
\lim_{\gb\to 0}\gl_{\gb}=\infty.
\end{equation}
\es
\Proof {\it Step 1: the monotonicity. }ÊLet $0<\gb_{1}<\gb_{2}$ and let $v_{\ge,1}$ and
 $v_{\ge,2}$ be  the corresponding solutions of (\ref{gen2}) with $\gb$ respectively replaced by $\gb_{1}$ and $\gb_{2}$. Since the $v_{\ge,i}$ are limit of solutions with finite boundary value there holds $v_{\ge,1}>v_{\ge,2}$
by comparison principle. Therefore
$$\gl_{\gb_{1}}:=\lim_{\ge\to 0}\ge v_{\ge,1}\geq \gl_{\gb_{2}}:=\lim_{\ge\to 0}\ge v_{\ge,2}.
$$ 
 Next, if we assume that there exist  $\gb_{i}$ ($i=1,2$) such that $0<\gb_{1}<\gb_{2}$ and $\gl_{\gb_{1}}= \gl_{\gb_{2}}=\gl$ and if $\gw_{1}$ and $\gw_{2}$ are the corresponding solutions of (\ref{clo3}) with $\gb=\gb_{i}$ and $\gl=\gl_{\gb_{1}}= \gl_{\gb_{2}}$, then (\ref{clo}) implies 
$$m^{-1}\gr(x)\leq\gw_{i}\leq m\gr(x)\quad\forall x\in S,
$$
for some $ m>0$. Set $\tilde\gw=\gw_{1}^{\gb_{2}/\gb_{1}}$, then
 \begin{equation}\label{clo4}\BA {l}
-div_{g}\left((\gb_{2}^2\tilde\gw^2+|\nabla \tilde\gw|^2)^{p/2-1}\nabla \tilde\gw\right)-\gb_{2}\gl(\gb_{2}^2\tilde\gw^2+|\nabla \tilde\gw|^2)^{p/2-1}\tilde\gw\\\phantom{-}
=(p-1)\left(1-\myfrac{\gb_{2}}{\gb_{1}}\right)
\left(\myfrac{\gb_{2}}{\gb_{1}}\right)^{p-1}\gw_{1}^{(p-1)(\gb_{2}/\gb_{1}-1)}\left(\gb_{1}^2\gw_{1}^2+|\nabla\gw_{1}|^2\right)^{(p-2)/2}\myfrac{|\nabla\gw_{1}|^2}{\gw_{1}}.
\EA\end{equation}
Therefore $\tilde\gw$ is a strict sub-solution. By homogeneity, and since $\prt_{\gn}\tilde\gw$ vanishes on $\prt S$, we can assume that $\tilde\gw\leq \gw_{2}$, that there exists $x_{0}\in S$ such that $\tilde\gw(x_{0})=\gw_{2}(x_{0})$ and the coincidence set of $\tilde\gw$ and $\gw_{2}$ is a subset of $S$. Let 
$$z=-\myfrac{1}{\gb_{2}}(\ln\gw_{2}-\ln\tilde\gw)=v_{2}-\tilde v.$$ 
Then $z\leq 0$, it is not identically zero, $z(x_{0})=0$ and $z(x)\to-\infty$ as $\gr(x)\to\prt S$. Because
$$\BA {l}-\Gd_{g}v_{2}-(p-2)\myfrac{D^2v_{2}\nabla v_{2}.\nabla v_{2}}{1+|\nabla v_{2}|^2}+\gb_{2}(p-1)|\nabla v_{2}|^2
\\[2mm]
\phantom{---------}
\leq
-\Gd_{g}\tilde v-(p-2)\myfrac{D^2\tilde v\nabla\tilde v.\nabla\tilde v}{1+|\nabla\tilde v|^2}+\gb_{2}(p-1)|\nabla\tilde v|^2
\EA$$
developing this inequality, we obtain that, at $x=x_{0}$, there holds
$$
\BA {l} \phantom{-----} -\Gd_{g}z-(p-2)\myfrac{D^2z\nabla v_{2}.\nabla v_{2}}
{1+|\nabla v_{2}|^2}
\\
\noalign{\medskip}
+(p-2)\left[\myfrac{D^2\tilde v\nabla\tilde v.\nabla\tilde v}{1+|\nabla\tilde v|^2}- \myfrac{D^2\tilde v\nabla v_2.\nabla v_2}{1+|\nabla v_2|^2}\right]
+\gb_{2}(p-1)\left[|\nabla v_2|^2-  |\nabla\tilde v|^2\right] \leq 0
\EA
$$
Since $\tilde v$, $v_2$ are $C^2$ in $S$,  the strong maximum principle yields  a contradiction. Therefore $\gb\mapsto\gl_{\gb}$ is decreasing.\smallskip

\noindent{\it Step 2: the continuity. }ÊLet $\{\gb_{n}\}$ be a positive sequence such that $\gb_{n}\to\gb_{0}$ and $v_{\gb_{n}}$ be the corresponding solution of
\begin{equation}\label{genn}\left\{ \BA {l}
 -\Gd_{g}v_{\gb_{n}}-(p-2)\myfrac{D^2v_{\gb_{n}} \nabla v_{\gb_{n}}.\nabla v_{\gb_{n}}}{1+|\nabla v_{\gb_{n}}|^2}+{\gb_{n}}(p-1)|\nabla v_{\gb_{n}}|^2=-\gl_{\gb_{n}}\quad\mbox{ in }S
 \\[2mm]\phantom{\left(\left(-----;-+v_{n} ----\right)^{p=x,}\nabla v_{0}\right)}
 \lim_{x\to\prt S}v_{\gb_{n}}( x)=\infty,
 \EA\right.
 \end {equation}
and let $v_{\ge,\gb_{n}}$ be the corresponding solutions of (\ref{gen2}) with $\gb=\gb_{n}$. Since $\ge v_{\ge,\gb_{n}}$ remains locally bounded in $S$ when $\gb_{n}$ remains in a compact subset of $(0,\infty)$ and converges to $\gl_{\gb_{n}}$ locally uniformly as $\ge\to 0$, the set $\{\gl_{\gb_{n}}\}$ is bounded. Up to a subsequence (not relabeled) we can assume that $\gl_{\gb_{n}}\to\bar\gl$ as $n\to \infty$. Thanks to \rife{clo} and \rife{clo2},  there holds
\begin{equation}\label{estn}
\abs {v_{\gb_{n}}+\myfrac{\ln\gr(x)}{\gb_{n}}}\leq C_0 \qquad\mbox{and }\;\abs {\nabla v_{\gb_{n}}}\leq \myfrac{C_1}{\gr(x)},
 \end {equation}
for some constants $C_0$, $C_1$, hence the sequence $\{v_{\gb_{n}}\}$ remains locally bounded in $W^{1,\infty}_{loc}(S)$ and, therfore, in $C^{2,\ga}_{loc}(S)$.  Up to a subsequence $v_{\gb_{n}}\to \bar v$ in $C^{2}_{loc}(S)$, and $\bar v$ is a solution of 
$$\left\{ \BA {l}
 -\Gd_{g}\bar v-(p-2)\myfrac{D^2\bar v \nabla \bar v.\nabla \bar v}{1+|\nabla \bar v|^2}+{\gb_{0}}(p-1)|\nabla \bar v|^2=-\bar\gl\quad\mbox{ in }S
 \\[2mm]\phantom{\left(\left(,-;-+v_{n} ----\right)^{p=x,}\nabla v_{0}\right)}
 \lim_{x\to\prt S}\bar v( x)=\infty.
 \EA\right.
$$
 By uniqueness of the  ergodic limit, $\bar\gl=\gl_{\gb_{0}}$, and  $\gl_{\gb_{n}}\to\gl_{\gb_{0}}$ for the whole sequence. \smallskip

\noindent{\it Step 3: (\ref{lim}) holds. } Let $\gw$ be a positive solution of 
 \begin{equation}\label{clo*}\left\{\BA {l}
-div_{g}\left((\gb^2\gw^2+|\nabla \gw|^2)^{p/2-1}\nabla \gw\right)=\gb\gl_{\gb}(\gb^2\gw^2+|\nabla \gw|^2)^{p/2-1}\gw\quad\mbox{in }S\\[2mm]
\phantom{;;-\left((\gb^2\gw^2+|\nabla \gw|^2)^{p/2-1}\nabla \gw_{0}\right)}
\gw=0\quad\mbox{on }\prt S.
\EA\right.\end{equation}
We normalize $\gw$ by
$$\myint{S}{}|\nabla\gw|^pdv_{g}=1.
$$
Therefore, if $\gm_{_{S}}$ is the first eigenvalue of $-div_{g}(|\nabla\,.\,|^{p-2}\nabla\,.)$ in $W_{0}^{1,p}(S)$, there holds
$$\myint{}{}|\gw|^pdv_{g}\leq \myfrac{1}{\gm_{_{S}}}.
$$
Multiplying (\ref{clo*}) by $\gw$ and integrating over $S$ yields to
\begin{equation}\label{dir1}
\myint{S}{}(\gb^2\gw^2+|\nabla\gw|^2)^{p/2}dv_{g}
=\gb(\gl_{\gb}+\gb)\myint{S}{}(\gb^2\gw^2+|\nabla\gw|^2)^{p/2-1}\gw^2dv_{g}.
\end{equation}
Clearly
$$\myint{S}{}(\gb^2\gw^2+|\nabla\gw|^2)^{p/2}dv_{g}\geq
\myint{S}{}|\nabla\gw|^{p}dv_{g}=1.
$$
If $p\geq 2$, 
$$\BA{l}\myint{S}{}(\gb^2\gw^2+|\nabla\gw|^2)^{p/2-1}\gw^2dv_{g}
\leq 2^{p/2-2}\myint{S}{}(\gw^{p}+\gw^2|\nabla\gw|^{p-2})dv_{g}\\\phantom{\myint{S}{}(\gb^2\gw^2+|\nabla\gw|^2)^{p/2-1}\gw^2dv_{g}}
\leq 2^{p/2-2}\left(1+\myfrac{2}{p}\right)\myint{S}{}\gw^{p}dv_{g}+2^{p/2-2}\left(1-\myfrac{2}{p}\right)\myint{S}{}|\nabla\gw|^{p})dv_{g}
\\\phantom{\myint{S}{}(\gb^2\gw^2+|\nabla\gw|^2)^{p/2-1}\gw^2dv_{g}}
\leq C_{_{p,S}}
\EA$$
This implies
\begin{equation}\label{p>2}
\gb(\gl_{\gb}+\gb)\geq \myfrac{1}{C_{_{p,S}}}\Longrightarrow
\gl_{\gb}\geq \myfrac{1}{C_{_{p,S}}\gb}-\gb.
\end{equation}
If $1<p<2$, 
$$\BA{l}\myint{S}{}\myfrac{\gw^2dv_{g}}{(\gb^2\gw^2+|\nabla\gw|^2)^{1-p/2}}
\leq \gb^{p-2}\myint{S}{}|\gw|^pdv_{g}\leq\myfrac{\gb^{p-2}}{\gm_{_{S}}}.
\EA$$
Therefore
\begin{equation}\label{p<2}
\gb^{p-1}(\gl_{\gb}+\gb)\geq \gm_{_{S}}\Longrightarrow
\gl_{\gb}\geq \myfrac{\gm_{_{S}}}{\gb^{p-1}}-\gb.
\end{equation}
Clearly (\ref{p>2}) and (\ref{p<2}) imply (\ref{lim}).\qeda\medskip
\vskip1em
\noindent{\Remark} Using the uniform ellipticity and the maximum principle,  (\ref{p>2}) and (\ref{p<2}) can possibly be improved in $
\gl_{\gb}\geq \myfrac{C}{\gb}$.
\vskip2em
We have now all the ingredients for the proof of Theorem B.
\vskip2em

\mbox {\bf Proof of Theorem B.}\quad If we set $\gw=e^{-\gb v}$ where $v$ is the solution of (\ref{gen1}), then $\gw$ is defined up to a multiplicative constant and satisfies  (\ref{clo*}). By \rlemma{regu}, $\gw\in C^1(\bar S)\cap C^2(S)$. Therefore the Theorem is obtained if we can prove that there exists a unique $\gb:=\gb_{_{S}}>0$ such that
\begin{equation}\label{O2}
\gl_{\gb}=\gb(p-1)+p-d-1.
\end{equation}
But the mapping $\gb\mapsto \gl_{\gb}-\gb(p-1)$ is continuous and decreasing on $(0,\infty)$. Clearly
$$\lim_{\gb\to\infty}\gl_{\gb}-\gb(p-1)=-\infty,
$$
and
$$\lim_{\gb\to0}\gl_{\gb}-\gb(p-1)=\infty,
$$
by \rprop{map}. The results follows by continuity.\qeda

\section{The regular case and Tolksdorf's result}
\setcounter{equation}{0}

If $\gb<0$, the equation satisfied by a a separable $p$-harmonic function $u$ under the form (\ref{s1}) is unchanged. However, if we set $\tilde\gb=-\gb$, then (\ref{s-p-h}) turns into
\begin{equation}\label{s-p-h*}
-div\left((\tilde\gb^2\gw+|\nabla'\gw|^2)^{p/2-1}\nabla'\gw\right)
=\tilde\gb(\tilde\gb(p-1)+N-p)(\tilde\gb^2\gw+|\nabla'\gw|^2)^{p/2-1}\gw.
\end{equation}
Furthermore, if a solution $\gw$ of (\ref{s-p-h*}) in $S\subset S^{N-1}$ exists which vanishes on $\prt  S$, then $\tilde\gb(p-1)+N-p>0$ by multiplying by $\gw$ and integration over $S$.
By setting 
$$v=-\myfrac{\ln\gw}{\tilde\gb},$$
then $v$ satisfies
 \begin{equation}\label{mainI*}\left\{ \BA {l}
 -div\left(\left(1+|\nabla' v|^2\right)^{p/2-1}\nabla' v\right)+\gb(p-1)\left(1+|\nabla' v|^2\right)^{p/2-1}|\nabla' v|^2\\[2mm]\phantom{---.\left(\left(1+|\nabla' v|^2\right)^{p/2-1}\right)}
 =-(\tilde\gb(p-1)+N-p)\left(1+|\nabla' v|^2\right)^{p/2-1}\quad\mbox{ in }S
 \\[2mm]\phantom{\left(\left(\right)^{p=x,2}\nabla' v\right)}
 \lim_{ \gs\to\prt S}v( \gs)=\infty.
 \EA\right.
 \end {equation}
 In the general setting of a Riemannian manifold, \rth{dev-th} and \rprop{map} are valid with $\gb$ replaced by $\tilde\gb$. The proof of  Theorem B holds except that (\ref{O2}) is replaced by
 \begin{equation}\label{O*2}
 \gl_{\tilde\gb}=\tilde\gb(p-1)+d+1-p.
 \end{equation}
 Because the function $\tilde\gb\mapsto \gl_{\tilde\gb}-\tilde\gb(p-1)$ is unchanged, the proof of Theorem B applies and shows that there exists a unique $\tilde\gb:=\tilde\gb_{_{S}}>0$ such that (\ref{O*2}) holds. Consequently we have proved the following result which contains Tolksdorf's initial result if $(M,g)=(S^{N-1},g_{0})$.
 
 \bcor{To*} Under the assumptions of \rth{dev-th} there exists a unique $\tilde\gb:=\tilde\gb_{_{S}}>0$ such that the problem
\begin{equation}\label{O1ter}\left\{\BA{l}
-div_{g}\left((\tilde\gb^2\gw^2+|\nabla \gw|^2)^{p/2-1}\nabla \gw\right)=\tilde\gb\left(\tilde\gb(p-1)+d+1-p\right)(\tilde\gb^2\gw^2+|\nabla \gw|^2)^{p/2-1}\gw\quad\mbox{in }S\\[2mm]
\phantom{;;-\left((\gb^2\gw_{0}^2+|\nabla \gw_{0}|^2)^{p/2-1}\nabla \gw_{0}\right)}
\gw=0\quad\mbox{on }\prt S,
\EA\right.
\end{equation}
admits a positive solution $\gw\in C^1(\bar S)\cap C^2(S)$. Furthermore $\gw$ is unique up to an homothethy.
\es

\appendix
\section{Appendix}
\setcounter{equation}{0}
We prove here the $C^{1,\gg}$ regularity up to the boundary, stated in \rlemma{regu}, for solutions of degenerate equations in divergence form
\begin{equation}\label{ap1}
\left\{\BA {l}
-div\left( a(x,u,\nabla u)\right)=B(x,u,\nabla u)\quad\mbox{in }S\\[2mm]
u=0\quad\mbox{on }\prt S.
\EA\right.
\end{equation}
We will assume that $a(x,s,\xi)$ satisfies the following conditions:
there exist constants $\la$, $\Lambda$,  $\beta >0$, and $\alpha \in (0,1]$, and a continuous function $\mu :S\times \BBR\to\BBR$ such that,   for every $s, t \in\BBR,$ for every $\xi, \eta \in \BBR^N$, and a.e. $x\in\Omega$:
\begin{equation}\label{tol1}
\frac{\partial a^i}{\partial \xi_j}(x,s,\xi)\,\eta_i\eta_j \geq \la (\mu(x,s)^2+|\xi|^2  )^{\frac{p-2}2}|\eta |^2,
\end{equation}
\begin{equation}\label{tol2}
\left| \frac{\partial a^i}{\partial \xi_j}(x,s,\xi)\right|  \leq \Lambda (\mu(x,s)^2+|\xi|^2  )^{\frac{p-2}2},
\end{equation}
\begin{equation}\label{xeu}
|a(x,s,\xi)-a(y,t,\xi)|\leq \beta \left(1+|\xi|^{p-2}+|\xi|^{p-1}  \right)  \,[ |x-y|^\alpha+|s-t|^{\alpha}],
\end{equation}
\begin{equation}\label{b}
|B(x,s,\xi)|\leq \beta (1+|\xi|^{p}).
\end{equation}
The model we have in mind is clearly
$$
a(x,u,\nabla u)= (\mu(x,u)^2+ |\nabla u|^2)^{\frac{p-2}2} \nabla u
$$ 
where $p>1$, and the function $\mu(x,s)$ is Lipschitz (or possibly H\"older) continuous. In many cases, as in the proof of \rlemma{regu}, the a priori information that $u$ is Lipschitz (or H\"older) continuous could allow us to consider only the case  $\mu=\mu(x)$.

The $C^{1,\gg}$ estimates, or similar kind of regularity results, are by now classical since the works of E. DiBenedetto \cite{Dib} and P. Tolksdorf \cite{To2} for the $p$--Laplace equation: as far as the global  regularity, up to the boundary, is concerned, we refer to the works of G. Lieberman (e.g. \cite{Lieb}) or to \cite{CDi}. Despite  a large amount of literature available, it seems that no exact reference applies to our model, so that, for the sake of completeness, we feel like giving  a proof of this result, at least detailing the possible slight modifications in order that previous results can be generalized. To this purpose, we observe that while the case $p\geq 2$ is somehow contained, if not in previous statements, at least  in previous arguments (specifically,  we refer to  \cite{Lieb}), this seems  not sure for the case $p<2$  because of
our growth assumption \rife{xeu} (roughly speaking, the $(x,s)$--derivatives may grow like $|\xi|^{p-2}$).  Finally, we note that the next result would still hold for a nonhomogeneous boundary condition ($u=\vfi$ on $\partial S$) provided $\vfi$ belongs to  $C^{1,\alpha}\left(\partial S \right)$.

\bth{liebe}
Let $S$ be a bounded $C^{1,\alpha}$ domain in $\BBR^N$, and assume that 
\rife{tol1}--\rife{b} hold true. If $u$ is  a bounded weak solution of \rife{ap1}, then there exists $\gg\in (0,1)$ 
%(depending on $\alpha$, \Lambda/\la,p,N)
such that 
$u\in C^{1,\gg}(\overline S)$ and moreover
$$
\|u\|_{C^{1,\gg}(\overline S)}\leq C\left(\Lambda/\la, \alpha, \|u\|_{\infty}, p,N,S  \right)\,.
$$
\es

\proof  Because our specific interest is in the boundary estimate, we only prove the regularity of $u$ around a point $x_0\in \partial S$ (the inner regularity is treated in the same manner). Up to straightening  the boundary, we can assume that locally $\partial S=\{x:x_N=0\}$ and $S=\{x: x_N>0\}$. 

We follow the standard approach via perturbation argument.
We denote   $B_R=\{x\,:\,|x-x_0|<R\}$,  $B_R^+=B_R\cap S$, and consider the solution $v$ of 
\begin{equation}\label{aux}
\left\{\BA {l}
-div\left( a(x_0,u(x_0),\nabla v)\right)=0\quad\mbox{in }B_R^+\\[2mm]
v=u\quad\mbox{on }\prt B_R^+.
\EA\right.
\end{equation}
Problem \rife{aux} has a unique solution $v\in W^{1,p}(B_R^+)$.
Due to assumptions \rife{tol1}--\rife{tol2}, the estimates concerning $v$ are well-established (\cite{Dib}, \cite{To2}, \cite{Lieb}). In particular, from Lemma 5 in \cite{Lieb} we have, for some $\sigma>0$,
\begin{equation}\label{oscv}
\mathop{{\rm osc}}\limits_{B_r^+} \nabla v \leq C\,\left( \frac rR\right)^{\sigma}\left(R^{-N}\int_{B_R^+}|\nabla v|^p  dx \right)^{\frac1p}\qquad\forall r<\frac R2
\end{equation}
where $C$, here and after,  depends only on the constants appearing in the hypotheses and possibly on $\|u\|_{\infty}$, in particular through the quantity $\sup \{|\mu(x,s)|\,,x\in \bar S, |s|\leq \|u\|_{\infty}\}$.  Moreover, since $a(x,s,\xi)\xi\geq c(|\xi|^p-|\mu|^p)$,  one easily deduces from \rife{aux}, using $v-u$ as test function  and Young's inequality, that
\begin{equation}\label{dv}
\int_{B_R^+} |\nabla v|^p\,dx\leq C \left( 1+ \int_{B_R^+} |\nabla u|^p\,dx\right)\,.
\end{equation}
Finally, the maximum  principle gives $\inf\limits_{\partial B_R^+} u\leq v\leq \sup\limits_{\partial B_R^+} u$, which yields
\begin{equation}\label{oscuv}
\mathop{{\rm osc}}\limits_{B_R^+} v \leq \mathop{{\rm osc}}\limits_{B_R^+} u\,.
\end{equation}
Now take $u-v$ as test function both in \rife{ap1} (restricted to $B_R^+$) and in \rife{aux} to obtain
$$
\int_{B_R^+} a(x,u,\nabla u)\nabla (u-v) \,dx - \int_{B_R^+}a(x_0,u(x_0),\nabla v)\nabla (u-v)\, dx = \int_{B_R^+} B(x,u,\nabla u)(u-v)\, dx\,.
$$
Denote $D_v\,:\,= \{x\in B_R^+\,:\, |\nabla u|< |\nabla v|\}$ and $D_u\,:\,=\{ x\in B_R^+\,:\, |\nabla v|\leq |\nabla u|\}$: hence we have
\begin{equation}\label{ap2}
\begin{array}{l}
\int_{D_v} \left[ a(x,u,\nabla u)-a(x,u,\nabla v)\right] \nabla (u-v) \,dx \\
\noalign{\medskip}
+ \int_{D_u}\left[ a(x_0,u(x_0), \nabla u)-a(x_0,u(x_0),\nabla v)\right] \nabla (u-v)\, dx \\
\noalign{\medskip}
\qquad = \int_{D_v} \left[ a(x_0,u(x_0),\nabla v)-a(x,u,\nabla v)\right] \nabla (u-v) \,dx\\
\noalign{\medskip}
+
\int_{D_u}\left[ a(x_0,u(x_0), \nabla u)-a(x,u,\nabla u)\right] \nabla (u-v)\, dx+ 
\int_{B_R^+} B(x,u,\nabla u)(u-v)\, dx
\end{array}
\end{equation}
Using \rife{xeu} and the definition of $D_v$, we have
$$
\begin{array}{l}
\int_{D_v} \left[ a(x_0,u(x_0),\nabla v)-a(x,u,\nabla v)\right] \nabla (u-v) \,dx\\
\noalign{\medskip}
\leq 2\beta \int_{D_v}  \left(1+|\nabla v|^{p-2}+|\nabla v|^{p-1}  \right) |\nabla v|  \,[ |x-x_0|^\alpha+|u(x)-u(x_0)|^{\alpha}]\,dx \\
\noalign{\medskip}
\leq C [ R^\alpha+(\mathop{{\rm osc}}\limits_{B_R^+} u )^\alpha]  
\int_{D_v}  \left(1+|\nabla v|^{p}  \right)  \,dx
\end{array}
$$
Similarly we estimate the second term in the right hand side of \rife{ap2}, and using also  \rife{b} we deduce
$$
\begin{array}{l}
\int_{D_v} \left[ a(x,u,\nabla u)-a(x,u,\nabla v)\right] \nabla (u-v) \,dx \\
\noalign{\medskip}
+ \int_{D_u}\left[ a(x_0,u(x_0), \nabla u)-a(x_0,u(x_0),\nabla v)\right] \nabla (u-v)\, dx \\
\noalign{\medskip}
\qquad \leq  
C [ R^\alpha+(\mathop{{\rm osc}}\limits_{B_R^+} u )^\alpha + \mathop{{\rm osc}}\limits_{B_R^+} u ]  
\int_{B_R^+}  \left(1+|\nabla v|^{p} +|\nabla u|^p \right)  \,dx\,,
\end{array}
$$
where we used that $\mathop{{\rm osc}}\limits_{B_R^+}( u-v) \leq 2\mathop{{\rm osc}}\limits_{B_R^+} u$ thanks to \rife{oscuv}.

Now, in both terms in the left hand side we use   \rife{tol1}  which implies, for every $(x,s,\xi)$:
\begin{equation}\label{coerc}
[a(x,s,\xi)-a(x,s,\eta)]\cdot(\xi-\eta)\geq c(\la)\left( \mu(x,s)^2+|\xi|^2+|\eta|^2\right)^{\frac{p-2}2}|\xi-\eta|^2\,.\end{equation}
If $p<2$ we get (recall that the generic constant $C$ may depend on $\|u\|_\infty$)
$$
\begin{array}{l}
\int_{D_v} \left[ a(x,u,\nabla u)-a(x,u,\nabla v)\right] \nabla (u-v) \,dx \\
\noalign{\medskip}
+ \int_{D_u}\left[ a(x_0,u(x_0), \nabla u)-a(x_0,u(x_0),\nabla v)\right] \nabla (u-v)\, dx \\
\noalign{\medskip}
\geq C \int_{D_v\cup D_u} \left[1+|\nabla u|^2+|\nabla v|^2\right]^{\frac{p-2}2} |\nabla (u-v) |^2\,dx 
\end{array}
$$
hence using H\"older inequality we end up with 
$$
\begin{array}{l}
\int_{B_R^+} | \nabla (u-v)|^p \,dx \leq  
C [ R^\alpha+(\mathop{{\rm osc}}\limits_{B_R^+} u )^\alpha + \mathop{{\rm osc}}\limits_{B_R^+} u ]^{q} 
\int_{B_R^+}  \left(1+|\nabla v|^{p} +|\nabla u|^p \right)  \,dx\,,
\end{array}
$$
with $q=\frac p2$.
If $p\geq 2$ we simply get rid of the term $\mu^2$ in \rife{coerc} and obtain the same inequality  with $q=1$.   Therefore, using also \rife{dv},  we conclude for any $p>1$
\begin{equation}\label{key}
\begin{array}{l}
\int_{B_R^+} | \nabla (u-v)|^p \,dx \leq  
C [ R^\alpha+(\mathop{{\rm osc}}\limits_{B_R^+} u )^\alpha + \mathop{{\rm osc}}\limits_{B_R^+} u ]^q 
\int_{B_R^+}  \left(1+ |\nabla u|^p \right)  \,dx
\end{array}
\end{equation}
with $q=\min(1,\frac p2)$.

Starting from the inequality \rife{key} it is possible to deduce the H\"older regularity of $\nabla u$ following well--known arguments. In particular, if  $u$ is Lipschitz continuous (as in our application in \rlemma{regu}) the conclusion is straightforward, since    \rife{key} implies
$$
\int_{B_R^+} | \nabla (u-v)|^p \,dx \leq  
C R^{N+ \alpha\, q}
$$ 
and   \rife{oscv}--\rife{dv} give that 
$\mathop{{\rm osc}}\limits_{B_r^+} \nabla v \leq C\,\left(   r/R\right)^{\sigma}$.

Then, defining  $(F)_r= \frac1{|B_r^+|}\int_{B_r^+}  F(y) dy$ for $F=\nabla u$ or $\nabla v$, we deduce 
$$
\begin{array}{l}
\int_{B_r^+} \left| \nabla u - \left(\nabla u\right)_r\right|^p\,dx \leq 
C \left[ \int_{B_r^+} \left|\nabla u- \nabla v\right|^p+  \int_{B_r^+} \left| \nabla v - \left(\nabla v\right)_r\right|^p\right]
\\ \noalign{\medskip}
\qquad\qquad\qquad \qquad
\leq  C[R^{N+\alpha q}+ r^N\, \left( \frac rR\right)^{\sigma\,p}]
\end{array}
$$
 and choosing $R=r^{\theta}$ for some  suitable $\theta<1$
 %($\theta=\frac{N+\sigma p}{N+\alpha q+\sigma p$),
%$$
%\int_{B_r^+} \left| \nabla u - \left(  \nabla u \right)_r\right|^p\,dx \leq Cr^{N+\gg\,p}
%$$
%for   $\gg=\sigma(1-\theta)\in (0,1)$. 
the conclusion follows from the results of Campanato \cite{Cam}.

In the general case, i.e.  when    a Lipschitz estimate on $u$ is not available, one need further work to estimate the right hand side of \rife{key}.  For this purpose, starting from \rife{key}, we can follow the arguments of   G. Lieberman (\cite{Lieb}, Section 3) and still get at the conclusion. 
\qed

%%%%%%%%%%%%%%%%%%%%%%%%%%
%%%%%%%%%%%%%%%%%%%%%%%%%%%%
%%BIBLIOGRAPHY%%%%%%%%%%%%%%%%%
%%%%%%%%%%%%%%%%%%%%%%%%%%%%


\begin{thebibliography}{777}

\bibitem{Berger} M. Berger, P. Gauduchon, E. Mazet, {\sl Le spectre d'une vari\'et\'e riemannienne.} Lecture Notes in Mathematics {\bf 194}. Springer--Verlag, Berlin/New York, 1971.

\bibitem{ber} S. Bernstein, \emph{Sur la g\'en\'eralization du probl\`eme de Dirichlet II}, Math. Ann. \textbf{69} (1910), 82-136.

\bibitem{BBoV1} M-F. Bidaut-V\'eron, R. Borghol, L. V\'eron,  \emph{Boundary Harnack inequality and a priori
estimates of singular solutions of quasilinear
elliptic equations}, Calc. Var. {\bf 27} (2006), 159-177.

\bibitem{BJV} M-F. Bidaut-V\'eron, M. Jazar, L. V\'eron,   \emph{Separable solutions of some quasilinear equations with source reaction}, J. Diff. Equ., to appear.

\bibitem{BoV} R. Borghol, L. V\'eron,  \emph{Boundary singularities of $N$-harmonic functions}, Comm. Part. Diff. Equ. {\bf 32} ( 2007), 1001-1015.

\bibitem{BoV2} R. Borghol, L. V\'eron,  \emph{Boundary singularities of solutions of $N$-harmonic equations with absorption}, J. Funct. Anal. {\bf  241} (2006), 611-637.

\bibitem{Cam} S. Campanato, \emph{Propriet\`a di H\"olderianit\`a di alcune classi di funzioni}, Annali Scuola Norm. Sup.   Pisa Cl. Sci. {\bf 17}, 175--188 (1963).

\bibitem{CDi} Chen Y. Z, E. DiBenedetto, \emph{Boundary estimates for solutions of nonlinear degenerate parabolic systems},  J. Reine Angew. Math. 395 (1989), 102-131.

\bibitem{Dib} E. DiBenedetto, \emph{$C\sp{1+\alpha }$ local regularity of weak solutions of degenerate elliptic equations}, Nonlinear Anal. {\bf 7}  (1983),  827--850. 

\bibitem{KV} S.  Kichenassamy and L. V\'eron. \emph{ Singular solutions of
the $p$-Laplace equation}, Math. Ann. {\bf 275} (1986), 599-615.

\bibitem{Kr} I. N.  Krol. \emph{The behaviour of the solutions of a certain quasilinear equation near 
zero cusps of the boundary}, Proc. Steklov Inst. Math.  {\bf 125} (1973), 130-136.

\bibitem{KM} I. N.  Krol and V. G. Mazja. \emph{The lack of
continuity and H\" older continuity of the solution of a certain
quasilinear equation}, Trudy Moskov. Mat. Ob\v s\v c. {\bf 26} (1972), 75-94.


\bibitem{doc} M.P. do Carmo {\sl Riemannian geometry},  Birkh\"auser Boston, Inc., Boston, MA, 1992.       

\bibitem{GT} D. Gilbarg , N. Trudinger, {\sl Partial Differential  Equations of Second Order}, 2nd ed., Springer-Verlag, Berlin/New-York, 1983.

\bibitem{pole} R.E. Greene, H. Wu  {\sl Function theory on manifolds which possess a pole}, Lecture Notes in Mathematics, 699. Springer, Berlin, 1979.

\bibitem{Ho} C.J. Holland,  \emph{ A new energy characterization of the smallest eigenvalue of the Schr\"odinger equation}, Comm. Pure Appl. Math. {\bf 30} (1977), 755-765. 


\bibitem{LL}  J.-M. Lasry,  P.-L. Lions,  \emph{  Nonlinear elliptic  equations with singular boundary conditions and stochastic control with state constraints. I.  The model problem},  Math. Ann. {\bf 283} (1989), 583-630.

\bibitem{Lieb} G. Lieberman, \emph{Boundary regularity for solutions of degenerate elliptic equations}, Nonlinear Anal. T.M.A. {\bf 12}, n.11 (1988), 1203--1219.

%\bibitem{L2}  P.-L. Lions,  \emph{ R\'esolution de probl\`emes elliptiques quasilin\'eaires}, Arch. Rat. Mech. An. {\bf 74} (1980), 234-254.

        
\bibitem {To} P. Tolksdorf,  \emph{On the Dirichlet problem for quasilinear equations in domains with conical boundary points}, Comm. Part. Diff. Equ. {\bf 8}  (1983), 773-817.

\bibitem {To2} P. Tolksdorf,  \emph{Regularity for a more general class of quasilinear elliptic equations}, J. Diff. Equ. {\bf 51}  (1984), 126-150.

%\bibitem{Se} J. Serrin, \emph{The problem of Dirichlet for quasilinear elliptic differential equations with many independent variables},
%Phil. Trans. Roy. Soc. London Ser. A {\bf 264}  (1969) 413-496. 

\bibitem{Ve1} L. V\'eron,  \emph{Some existence and uniqueness results for solutions of some quasilinear elliptic equations on compact Riemannian manifolds}, Colloquia Mathematica Societatis J\'anos Bolyai {\bf 62}  (1991), 317-352.
\end{thebibliography}
\end{document}